\documentclass[11pt]{article}
\usepackage{graphicx}
\usepackage{latexsym,amssymb}
\input epsf

\title{\textbf{Edgewise Cohen-Macaulay connectivity of
partially ordered sets}}

\author{Christos~A.~Athanasiadis\\
\small Department of Mathematics\\
\small University of Athens\\
\small Panepistimioupolis \\
\small Athens 15784, Hellas (Greece)\\
\small \texttt{caath@math.uoa.gr}
\vspace{5pt}
\and
Myrto~Kallipoliti\\
\small Department of Mathematics\\
\small University of Vienna\\
\small Oskar-Morgenstern-Platz 1\\
\small 1090 Vienna, Austria\\
\small \texttt{myrto.kallipoliti@univie.ac.at} }

\date{\small November 6, 2015}

  \def\kk{{\mathbf k}}
  
  \def\bB{{\mathcal B}}
  \def\iI{{\mathcal I}}
  \def\jJ{{\mathcal J}}
  \def\eE{{\mathcal E}}
  \def\lL{{\mathcal L}}
  
  \def\pP{{\mathcal P}}
  \def\qQ{{\mathcal Q}}
  
  \def\rR{{\mathcal R}}
  \def\ast{{\rm st}}
  \def\cost{{\rm cost}}

  \def\lk{{\rm link}}
  \def\rank{{\rm rank}}
  
  \def\sm{\smallsetminus}

  \newcommand{\qed}{$\hfill \Box$}

\begin{document}
\maketitle

\newtheorem{theorem}{Theorem}[section]
\newtheorem{proposition}[theorem]{Proposition}
\newtheorem{corollary}[theorem]{Corollary}
\newtheorem{definition}[theorem]{Definition}
\newtheorem{remark}[theorem]{Remark}
\newtheorem{lemma}[theorem]{Lemma}
\newtheorem{example}[theorem]{Example}
\newtheorem{examples}[theorem]{Examples}
\newtheorem{conjecture}[theorem]{Conjecture}
\newtheorem{fact}[theorem]{Fact}
\newtheorem{question}[theorem]{Question}
\newtheorem{observation}[theorem]{Observation}
\newtheorem{claim}[theorem]{Claim}

\begin{abstract}
The proper parts of face lattices of convex polytopes are shown to satisfy
a strong form of the Cohen--Macaulay property, namely that removing
from their Hasse diagram all edges in any closed interval results in
a Cohen--Macaulay poset of the same rank. A corresponding notion of
edgewise Cohen--Macaulay connectivity for partially ordered sets is
investigated. Examples and open questions are discussed.

\end{abstract}


\section{Introduction and results}
\label{sec:intro}

Cohen--Macaulay simplicial complexes and partially ordered sets (posets)
have been studied intensely for the past few decades, mainly due to their
importance in several areas of mathematics; see~\cite{Bj14} for an informal
discussion and references.

Cohen--Macaulay connectivity, introduced by Baclawski~\cite{Ba82}, is a
sophisticated notion of connectivity for simplicial complexes and posets;
it includes as a special case the notion of vertex-connectivity of graphs.
This paper investigates a new notion of connectivity for posets which
is motivated by that of edge-connectivity of graphs. Given a positive
integer $k$, a poset $\pP$ is called \emph{edgewise $k$-Cohen--Macaulay}
if $\pP$ is Cohen--Macaulay of rank at least $k-1$ and the removal from
the Hasse diagram of $\pP$ of all edges in any closed interval of $\pP$ of
rank less than $k$ results in a Cohen--Macaulay poset of the same rank as
$\pP$. The main problem we propose to study asks to determine the largest
integer $k$ for which a given poset is edgewise $k$-Cohen--Macaulay; this
integer is called the edgewise Cohen--Macaulay connectivity.

The main result of this paper solves this problem for one important class
of posets, namely the proper parts of Gorenstein* lattices. We call
\emph{edgewise strongly Cohen--Macaulay} the posets whose
edgewise Cohen--Macaulay connectivity is maximum possible, equal to one
more than their rank. Thus, edgewise strongly Cohen--Macaulay posets have
the property that removing from their Hasse diagram all edges in any closed
interval results in a Cohen--Macaulay poset of the same rank. Recall that
the class of Gorenstein* lattices includes all face lattices of convex
polytopes and that the proper part $\overline{\lL}$ of a lattice $\lL$ is
the poset obtained from $\lL$ by removing its minimum and maximum elements.

\begin{theorem} \label{thm:goren}
The poset $\overline{\lL}$ is edgewise strongly Cohen--Macaulay for every
Gorenstein* lattice $\lL$.
\end{theorem}

This theorem is proven in Section~\ref{sec:goren} using methods of
topological combinatorics. The notion of edgewise Cohen--Macaulay
connectivity, adopted in this paper, is introduced in
Section~\ref{sec:ecm}, where elementary properties of edgewise
$k$-Cohen--Macaulay posets are established, examples are listed and a
comparison to $k$-Cohen--Macaulay posets, in the sense of
Baclawski~\cite{Ba82}, is given for $k=2$
(Proposition~\ref{prop:2CM2ECM}). Basic definitions and notation
needed to understand this paper are summarized in Section~\ref{sec:pre}.
Some open problems and related remarks are discussed in
Section~\ref{sec:rem}.

\section{Preliminaries}
\label{sec:pre}

This section fixes notation and recalls basic definitions and background
on the combinatorics and topology of simplicial complexes and partially
ordered sets. For more information on these topics and any undefined terms,
the reader is referred to the sources~\cite{Bj95} \cite{StaCCA}
\cite[Chapter~3]{StaEC1} \cite{Wa07}. Basic background on algebraic
topology can be found in~\cite{Mu84}.

\medskip
\noindent
\textbf{Simplicial complexes}.
An (abstract) \emph{simplicial complex} on a ground set $E$ is
a collection $\Delta$ of subsets of $E$, called \emph{faces}, such that
$F \subseteq G \in \Delta$ implies $F \in \Delta$. The elements $v \in E$
for which $\{v\} \in \Delta$ are the \emph{vertices} of $\Delta$. The dimension of
a face $F \in \Delta$ is defined as one less than the cardinality of $F$.
The dimension of $\Delta$, denoted $\dim(\Delta)$, is the maximum dimension
of its faces. This paper is concerned with order complexes of finite posets
(discussed later in this section); in particular, all simplicial complexes
considered here will be finite.

Every simplicial complex $\Delta$ has a geometric
realization~\cite[Section~9]{Bj95}, uniquely defined up to homeomorphism;
it will be denoted by $|\Delta|$. All topological properties of $\Delta$
mentioned in the sequel refer to those of $|\Delta|$.

A subcollection $\Gamma \subseteq \Delta$ which is a simplicial complex
on its own is called a \emph{subcomplex}. The geometric realization
$|\Gamma|$ of such a complex may be considered as a subspace of $|\Delta|$.
The \emph{link} of a face $F \in \Delta$ is the subcomplex defined as
$\lk_\Delta (F) = \{ G \sm F: F \subseteq G \in \Delta \}$. The \emph{open
star} of a face $F \in \Delta$ is defined as $\ast_\Delta(F) = \{G \in \Delta: F
\subseteq G\}$. Its complement in $\Delta$, denoted $\cost_\Delta(F)$, is
a subcomplex of $\Delta$ called the \emph{contrastar} of $F$; it consists
of all faces of $\Delta$ which do not contain $F$. A subcomplex $\Gamma$
of $\Delta$ is said to be \emph{vertex-induced} if $\Gamma$ contains all
faces $F \in \Delta$ for which all elements of $F$ are vertices of $\Gamma$.
Given $A \subseteq E$, the (induced) subcomplex consisting of all faces
of $\Delta$ which are disjoint from $A$ will be denoted by $\Delta \sm A$.
The complex $\Delta$ is said to be a \emph{cone} with apex $v \in E$ if
$F \cup \{v\} \in \Delta$ for every $F \in \Delta$; such a complex is
contractible.

The following standard lemma will be used in Section~\ref{sec:goren}.

\begin{lemma} {\rm (\cite[Lemma~70.1]{Mu84})} \label{lem:induced}
Let $\Delta$ be a simplicial complex and $\Gamma$ be a vertex-induced
subcomplex with vertex set $A$. Then the space $|\Delta \sm A|$ is a
deformation retract of $|\Delta| \sm |\Gamma|$.
\end{lemma}

A simplicial complex $\Delta$ is \emph{Cohen--Macaulay} over a field
$\kk$ if $\widetilde{H}_i \, (\lk_\Delta(F), \kk) = 0$ for all $F \in
\Delta$ (including the empty face) and $i < \dim \, (\lk_\Delta(F))$,
where $\widetilde{H}_* (\Gamma, \kk)$ denotes reduced simplicial
homology of $\Gamma$ with coefficients in $\kk$, and \emph{Gorenstein*}
over $\kk$ if it is Cohen--Macaulay over $\kk$ and $\widetilde{H}_i \,
(\lk_\Delta(F), \kk) = \kk$ for all $F \in \Delta$ and $i = \dim \,
(\lk_\Delta(F))$. Given a positive integer $k$, the complex $\Delta$ is
\emph{$k$-Cohen-Macaulay} (or \emph{doubly Cohen-Macaulay}, for $k=2$)
over $\kk$ if $\Delta \sm A$ is Cohen--Macaulay over $\kk$ of the same
dimension as $\Delta$ for every subset $A$ of the ground set of
cardinality less than $k$ (including $A = \varnothing$).

Every Cohen--Macaulay simplicial complex $\Delta$ is \emph{pure}, meaning
that all maximal (with respect to inclusion) faces of $\Delta$ have
dimension equal to $\dim(\Delta)$. The properties of being
Cohen--Macaulay, Gorenstein* or doubly Cohen--Macaulay are all topological,
i.e., they depend only on the homeomorphism type of $|\Delta|$ and the
field $\kk$; more precise statements appear in~\cite{StaCCA, Wa81}.
From now on, we consider the field $\kk$ fixed and suppress it
from our terminology. Our results continue to hold if $\kk$ is replaced
by the ring of integers.

\begin{figure}
  \epsfysize = 1.2 in \centerline{\epsffile{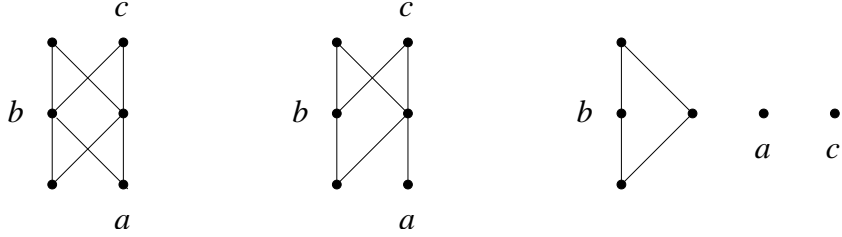}}
  \caption{A poset $\pP$ and the posets $\pP \ominus [a, b]$ and $\pP
  \ominus [a, c]$}
  \label{fig:ecmA}
  \end{figure}

\medskip
\noindent
\textbf{Partially ordered sets}. This paper considers only finite
partially ordered sets (posets). Thus, a poset is a pair $\pP = (P,
\le)$ consisting of a finite set $P$ and a partial order (reflexive,
antisymmetric and transitive binary relation) $\le$ on $P$. We will
denote by $\eE(\pP)$ the set of all \emph{cover relations} of $\pP$,
meaning ordered pairs $(a, b) \in P \times P$ with $a < b$ for
which no $x \in P$ satisfies $a < x < b$. The \emph{Hasse diagram}
of $\pP$ is the directed graph with vertex set $P$ and edge set
$\eE(\pP)$. We assume familiarity with drawing posets via their
Hasse diagrams; see Figure~\ref{fig:ecmA} on the left for an example
of a poset with six elements and eight cover relations.

Given elements $a, b \in P$ with $a \le b$, the closed interval
$[a, b]_\pP$ is defined as the (induced) subposet of $\pP$ on the set
$\{x \in P: a \le x \le b\}$, where the subscript $\pP$ can be dropped
if there is no danger of confusion. Similar remarks apply to open
and half-open intervals. The subposets on the sets $\{x \in
P: x \ge a\}$ and $\{x \in P: x \le b\}$ will be denoted by $\pP_{\ge
a}$ and $\pP_{\le b}$, respectively. The set of two-element closed
intervals (in other words, the set of undirected edges of the Hasse
diagram) of $\pP$ will be denoted by $E(\pP)$.

A \emph{chain} (respectively, \emph{antichain}) in $\pP$ is any
subset of $P$ which consists of pairwise comparable (respectively,
incomparable) elements. The poset $\pP$ is \emph{graded} if all
maximal (with respect to inclusion) chains in $\pP$ have the same
length (defined as one less than their cardinality). This common
length is then called the \emph{rank} of $\pP$ and denoted by
$\rank(P)$.

The poset $\pP$ is \emph{bounded} if it has a minimum element
$\hat{0}$ and a maximum element $\hat{1}$; the elements of $P$
covering $\hat{0}$ are called \emph{atoms} and those covered by
$\hat{1}$ are called \emph{coatoms}. The \emph{proper part} of a
bounded poset $\pP$ is defined as the subposet $\overline{\pP}$
obtained from $\pP$ by removing $\hat{0}$ and $\hat{1}$. Conversely,
the poset obtained from a poset $\qQ$ by artificially adding a
minimum element $\hat{0}$ and a maximum element $\hat{1}$
is denoted by $\widehat{\qQ}$. The \emph{ordinal sum} of two
posets $(P, \le_P)$ and $(Q, \le_Q)$ on disjoint ground sets $P$
and $Q$ is the set $P \cup Q$ partially ordered by setting $x \le
y$ if $x, y \in P$ and $x \le_P y$, or $x, y \in Q$ and $x \le_Q y$,
or $x \in P$ and $y \in Q$, for $x, y \in P \cup Q$.

A \emph{lattice} is a poset any two elements $x, y$ of which have a
least upper bound, called the \emph{join} and denoted by $x \vee y$,
and a greatest lower bound, called the \emph{meet} and denoted by
$x\wedge y$. Cleary, every finite lattice is bounded. Such a lattice
$\lL = (L, \le)$ is called \emph{atomic} if every element of $L$ is
the join of atoms, \emph{relatively atomic} if every closed interval
of $\lL$ is atomic and \emph{semimodular} if for all $x, y \in L$,
the join $x\vee y$ covers $y$ whenever $x$ covers the meet $x\wedge
y$. A lattice which is both atomic and semimodular is called
\emph{geometric}.

The \emph{order complex} associated to a poset $\pP = (P,\le)$ is
the simplicial complex on the ground set $P$ whose faces
are the chains of $\pP$; it is denoted $\Delta(\pP)$. Note that $\pP$
is graded if and only if $\Delta(\pP)$ is pure; in that case, $\rank
(\pP) = \dim(\Delta(\pP))$. The poset $\pP$ is called Cohen--Macaulay,
Gorenstein* or $k$-Cohen--Macaulay if $\Delta(\pP)$ has
the corresponding property. An equivalent definition for $\pP$ to be
Cohen--Macaulay is that $\widetilde{H}_i \, (\Delta(\jJ), \kk) = 0$
for all open intervals $\jJ$ in $\widehat{\pP}$ and indices $i <
\rank(\jJ)$. A lattice $\lL$ is said to be Gorenstein* (by slight
abuse of language) if its proper part $\overline{\lL}$ is a
Gorenstein* poset. The class of Gorenstein* lattices includes all
face lattices of convex polytopes and, more generally, regular cell
decompositions of the sphere having the intersection property.

\section{Edgewise Cohen--Macaulay connectivity}
\label{sec:ecm}

This section introduces the notion of edgewise Cohen-Macaulay
connectivity for posets, establishes some of its elementary
properties and compares it with the standard notion of Cohen-Macaulay
connectivity, due to Baclawski~\cite{Ba82}.

Given a poset $\pP = (P, \le)$ and a closed interval $\iI \subseteq
\pP$, we denote by $\pP \ominus \iI$ the poset whose Hasse diagram
is obtained from that of $\pP$ by removing all edges within $\iI$.
Equivalently, $\pP \ominus \iI$ is the unique poset on the ground
set $P$ which satisfies $\eE(\pP \ominus \iI) = \eE(\pP) \sm \eE(\iI)$.
An example is shown in Figure~\ref{fig:ecmA}. Note that the elements of
$\iI$ are pairwise incomparable in $\pP \ominus \iI$.

The main definition of this paper is as follows.
\begin{definition} \label{def:kECM}
Let $k$ be a positive integer. A Cohen--Macaulay poset $\pP$ is called
\begin{itemize}
\itemsep=0pt
\item[{\rm (a)}]
edgewise $k$-Cohen--Macaulay if $\rank(P) \ge k-1$ and $\pP \ominus \iI$
is Cohen--Macaulay of the same rank as $\pP$ for every closed interval
$\iI \subseteq \pP$ of rank less than $k$; and

\item[{\rm (b)}]
edgewise strongly Cohen--Macaulay if $\pP \ominus \iI$ is Cohen--Macaulay
of the same rank as $\pP$ for every closed interval $\iI \subseteq \pP$.
\end{itemize}
The largest integer $k$ for which $\pP$ is edgewise $k$-Cohen--Macaulay
is called the edgewise Cohen-Macaulay connectivity of $\pP$.
\end{definition}

In particular, the edgewise 1-Cohen--Macaulay posets are exactly the
Cohen--Macaulay posets, whereas a poset $\pP$ is edgewise 2-Cohen--Macaulay
(or edgewise doubly Cohen--Macaulay) if it is Cohen--Macaulay of positive
rank and $\pP \ominus e$ is Cohen--Macaulay of the same rank as $\pP$ for
every $e \in E(\pP)$.

\begin{example} \label{ex:basic} \rm
(a) A poset of rank one is Cohen--Macaulay (respectively, 2-Cohen--Macaulay
or edgewise 2-Cohen--Macaulay) if and only if it is connected (respectively,
2-vertex-connected or 2-edge-connected) as a graph.

(b) Chains are Cohen--Macaulay but neither 2-Cohen--Macaulay, nor edgewise
2-Cohen--Macaulay; their edgewise Cohen-Macaulay connectivity is equal to one.

(c) The poset on the left of Figure~\ref{fig:ecmA} is 2-Cohen--Macaulay and
edgewise 2-Cohen--Macaulay but neither 3-Cohen--Macaulay, nor edgewise
3-Cohen--Macaulay. Thus, its edgewise Cohen-Macaulay connectivity is equal
to two; see Remark~\ref{rem:ordinal} for a generalization.

(d) Let $\bB_n$ be the Boolean lattice of rank $n$
\cite[Example~3.1.1]{StaEC1}, consisting of all subsets of $\{1, 2,\dots,n\}$,
ordered by inclusion. It follows from Theorem~\ref{thm:goren} that $\overline
{\bB}_n$ is edgewise strongly Cohen--Macaulay for every $n$.
\qed
\end{example}

The following elementary properties of edgewise $k$-Cohen--Macaulay posets
are similar to corresponding properties of $k$-Cohen--Macaulay posets.

\begin{proposition} \label{prop:elementary}
\begin{itemize}
\itemsep=0pt
\item[{\rm (a)}]
The class of edgewise $(k+1)$-Cohen--Macaulay posets is strictly
included in that of edgewise $k$-Cohen--Macaulay posets for all
$k \ge 1$.

\item[{\rm (b)}]
Every poset which is edgewise $k$-Cohen--Macaulay for some $k \ge
2$ contains at least two minimal and at least two maximal elements.

\item[{\rm (c)}]
Let $\pP$ be edgewise $k$-Cohen--Macaulay. Then every open interval
$\jJ$ in $\widehat{\pP}$ is edgewise $r$-Cohen--Macaulay for $r =
\min\{k, \rank(\jJ) + 1\}$.
\end{itemize}
\end{proposition}

\noindent
\emph{Proof.} Part (a) is a direct consequence of
Definition~\ref{def:kECM}. The inclusion is strict since, by the same
definition, edgewise $k$-Cohen--Macaulay posets of rank $k-1$ are not
edgewise $(k+1)$-Cohen--Macaulay.

To verify part (b), suppose $\pP$ is Cohen--Macaulay with a
unique minimal (hence minimum) element $\hat{0}$, or a unique maximal
(hence maximum) element $\hat{1}$ and note that removing any edge
which contains $\hat{0}$ or $\hat{1}$ from the Hasse diagram of $\pP$
results in a poset which is either not graded, or which has rank less
than that of $\pP$. Thus, such a poset is not edgewise doubly
Cohen--Macaulay and hence not edgewise $k$-Cohen--Macaulay for any
$k \ge 2$ either.

To prove part (c) note that $\jJ$ is Cohen--Macaulay, as an interval
in a Cohen--Macaulay poset, and consider any closed interval $\iI$ of
rank less than $k$ in $\jJ$. Then $\iI$ is a closed interval of rank
less than $k$ in $\pP$ as well and thus $\pP \ominus \iI$ is
Cohen--Macaulay of the same rank as $\pP$. As a result, $\widehat{\pP}
\ominus \iI$ is Cohen--Macaulay of the same rank as $\widehat{\pP}$.
Being an open interval in $\widehat{\pP} \ominus \iI$, the poset $\jJ
\ominus \iI$ must be Cohen--Macaulay of the same rank as $\jJ$ and the
proof follows.
\qed

\bigskip
The next statement exhibits a large class of edgewise doubly
Cohen--Macaulay posets; it implies, for instance, that $\pP$ is edgewise
doubly Cohen--Macaulay whenever $|\Delta(\pP)|$ is homeomorphic to a
sphere of positive dimension.

\begin{proposition} \label{prop:2CM2ECM}
Every doubly Cohen--Macaulay poset of positive rank is edgewise doubly
Cohen--Macaulay.
\end{proposition}

To prepare for the proof of Proposition~\ref{prop:2CM2ECM}, we make
the following observation. Recall from Section~\ref{sec:pre} that
$\cost_{\Delta(\pP)} (F)$ is the simplicial complex of chains in a
poset $\pP$ which do not contain a given chain $F \in \Delta(\pP)$.
Since no chain in $\pP \ominus [a, b]$ can contain two distinct
elements of $[a, b]$, we have
\begin{equation} \label{eq:costinclusion}
\Delta(\pP \ominus [a, b]) \ \subseteq \bigcap_{a \le x < y \le b}
\cost_{\Delta(\pP)} \, (\{x, y\})
\end{equation}
for every closed interval $[a, b] \subseteq \pP$.

\begin{lemma} \label{lem:cost}
Let $\pP$ be a finite poset which contains no three-element closed
interval. Then $\Delta(\pP \ominus e) = \cost_{\Delta(\pP)} (e)$ for
every $e \in E(\pP)$.
\end{lemma}

\noindent
\emph{Proof.} Let $e \in E(\pP)$. We claim that for every pair $\{u, v\}$
of elements of $\pP$ other than $e$ with $u < v$, there exists a saturated
chain in $\pP$ with minimum element $u$ and maximum element $v$ which does
not contain $e$. Indeed, this is trivial if $\{u, v\} \in E(\pP)$.
Otherwise, one can consider any maximal chain of the closed interval
$[u, v]$ of $\pP$ and in case this chain contains both elements of $e$,
since $\pP$ contains no three-element closed intervals, one can replace
one of the elements of $e$ with other elements of $\pP$ to obtain a chain
with the desired properties.

As a result of the claim, if $u < v$ holds in $\pP$ for some pair $(u, v)$
of elements of $\pP$ other than $e$, then $u < v$ holds in $\pP \ominus e$
as well. This implies that every chain in $\pP$ which does not contain
$e$ is also a chain in $\pP \ominus e$, in other words that
$\cost_{\Delta(\pP)} (e) \subseteq \Delta(\pP \ominus e)$. The reverse
inclusion is a special case of~(\ref{eq:costinclusion}).
\qed

\bigskip
\noindent
\emph{Proof of Proposition~\ref{prop:2CM2ECM}.} Let $\pP$ be a doubly
Cohen-Macaulay poset and let $e \in E(\pP)$. We need to show that
$\Delta(\pP \ominus e)$ is Cohen--Macaulay of the same dimension as
$\Delta(\pP)$. Since the class of doubly Cohen--Macaulay complexes is
closed under taking links of faces, no open interval in $\pP$ can be a
singleton. Thus, from Lemma~\ref{lem:cost} we get $\Delta(\pP \ominus
e) = \cost_{\Delta(\pP)} (e)$. Moreover, since $\Delta(\pP)$ is doubly
Cohen--Macaulay, \cite[Proposition~2.8]{MM14} and its proof imply
that $\cost_{\Delta(\pP)} (e)$ is Cohen--Macaulay of the same dimension
as $\Delta(\pP)$ and the proof follows.
\qed

  \begin{figure}
  \epsfysize = 0.65 in \centerline{\epsffile{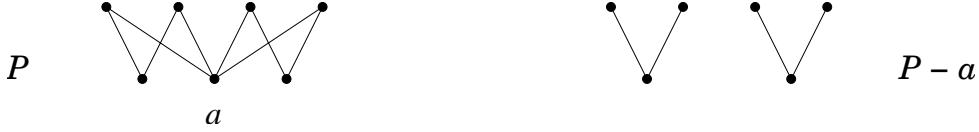}}
  \caption{An edgewise 2-Cohen--Macaulay poset $\pP$ which is not
  2-Cohen--Macaulay}
  \label{fig:ecmB}
  \end{figure}

\medskip
We conclude this section with the following remarks.

\begin{figure}
  \epsfysize = 1.15 in \centerline{\epsffile{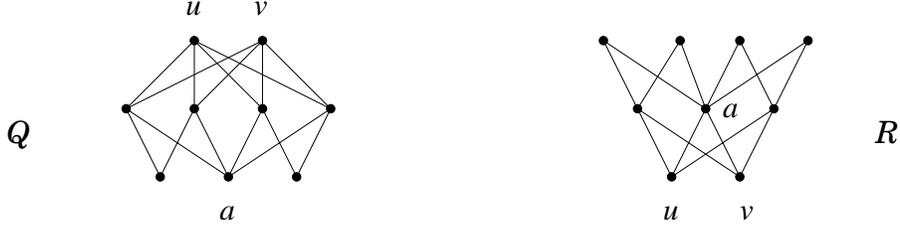}}
  \caption{Two homeomorphic posets}
  \label{fig:ecmC}
  \end{figure}

\begin{remark} \label{rem:converse} \rm
(a) The converse of Proposition~\ref{prop:2CM2ECM} is false.
Figure~\ref{fig:ecmB} shows an edgewise 2-Cohen--Macaulay poset $\pP$
of rank one which is not 2-Cohen--Macaulay. For an example of rank two,
consider the ordinal sum $\qQ$ of $\pP$ with a two-element antichain,
shown on the left of Figure~\ref{fig:ecmC}. One can easily check that
$\qQ \ominus e$ is pure-shellable of rank two for every $e \in E(\qQ)$
and conclude that $\qQ$ is edgewise 2-Cohen--Macaulay. On the other
hand, $\qQ \sm a$ is not Cohen--Macaulay since, for instance, removing
its two maximal elements results in a rank-selected subposet which is
disconnected of rank one and hence not Cohen--Macaulay. Therefore,
$\qQ$ is not 2-Cohen--Macaulay.

We leave it to the interested reader to check that, more generally,
the ordinal sum of $\pP$ and any number of two-element antichains,
taken in this order, is an edgewise 2-Cohen--Macaulay poset which is
not 2-Cohen--Macaulay.

(b) The converse of Proposition~\ref{prop:2CM2ECM} does hold for the
class of semimodular lattices. For a semimodular lattice $\lL$, the
following conditions are in fact equivalent: (i) $\lL$ is geometric; (ii)
$\overline{\lL}$ is 2-Cohen--Macaulay; and (iii) $\overline{\lL}$ is
edgewise 2-Cohen--Macaulay. The equivalence ${\rm (i)} \Leftrightarrow
{\rm (ii)}$ is the content of \cite[Theorem~3.1]{Ba82}. Suppose that
$\lL$ is semimodular but not geometric. Then, as discussed in the proof
of \cite[Theorem~3.1]{Ba82}, there exists an element $b \in \lL$ which
is not an atom and covers a unique element $a$. As a result, either
$\overline{\lL}$ has a unique maximal element (equal to $a$, if $b =
\hat{1}$), or else $\lL$ has an open interval of rank one which
contains $a$ and $b$ and hence is not edgewise 2-Cohen--Macaulay.
Proposition~\ref{prop:elementary} implies that $\overline{\lL}$ is not
edgewise 2-Cohen--Macaulay either. We have shown that ${\rm (iii)}
\Rightarrow {\rm (i)}$. The implication ${\rm (ii)} \Rightarrow
{\rm (iii)}$ follows from Proposition~\ref{prop:2CM2ECM}.
\end{remark}

\begin{remark} \label{rem:not-topo} \rm
(a) Contrary to the situation with (double) Cohen--Macaulayness,
edgewise double Cohen--Macaulayness is not a topological property.
Consider, for instance, the posets $\qQ$ and $\rR$
shown in Figure~\ref{fig:ecmC}. Both order complexes $\Delta(\qQ)$ and
$\Delta(\rR)$ are homeomorphic to the suspension over the wedge of two
circles. However, $\qQ$ is edgewise 2-Cohen--Macaulay, as already
commented in Remark~\ref{rem:converse} (a), and $\rR$ is not, since the
poset $(\rR \ominus [u, a])_{>u}$ is disconnected of rank one and hence
$\rR \ominus [u, a]$ is not Cohen--Macaulay.

Consideration of the ordinal sum of $n-1$ two-element antichains and
of the proper part $\overline{\bB}_n$ of the Boolean lattice of rank $n$
shows that edgewise $k$-Cohen--Macaulayness for $k \ge 3$ and edgewise
strong Cohen--Macaulayness are not topological properties either.

(b) The ordinal sum of edgewise doubly Cohen--Macaulay (respectively,
edgewise strongly Cohen--Macaulay) posets is not always edgewise doubly
Cohen--Macaulay (respectively, edgewise strongly Cohen--Macaulay). For
a counterexample, take the ordinal sum of two two-element antichains
and the poset of Figure~\ref{fig:ecmB} and argue as with the poset
$\rR$ of part (a).
\end{remark}

\begin{remark} \label{rem:ordinal} \rm
The natural generalizion of Proposition~\ref{prop:2CM2ECM} to
$k$-Cohen--Macaulayness fails for $k \ge 3$. Indeed, let $\pP$ be the
ordinal sum of $n$ copies of a $k$-element antichain. This poset
(shown on the left of Figure~\ref{fig:ecmA} for $n=3$ and $k=2$)
is $k$-Cohen--Macaulay, since removing fewer than $k$ elements from
$\pP$ yields the ordinal sum of $n$ nonempty antichains, which is
Cohen--Macaulay, as an ordinal sum of Cohen--Macaulay posets, of the
same rank as $\pP$. However, $\pP$ is not edgewise $k$-Cohen--Macaulay
for $k \ge 3$ since $\pP\ominus \iI$ has a disconnected Hasse diagram
for every rank two closed interval $\iI$ which contains a maximal or
minimal element of $\pP$.
\end{remark}

\section{Gorenstein* lattices}
\label{sec:goren}

This section proves Theorem~\ref{thm:goren}. The following lemma
provides a large class of posets $\pP$ for which equality holds in
(\ref{eq:costinclusion}) for every interval $[a, b] \subseteq \pP$;
it applies to geometric and Gorenstein* lattices (the latter have a
nowhere zero M\"obius function and are therefore relatively atomic,
for instance, by \cite[Corollary~3.9.5]{StaEC1}).

\begin{lemma} \label{lem:goren}
Let $\lL$ be a graded, relatively atomic lattice. Then
\begin{equation} \label{eq:lematomic}
\Delta(\overline{\lL} \ominus [a, b]) \ = \bigcap_{a \le x < y \le b}
\cost_{\Delta(\overline{\lL})} \, (\{x, y\})
\end{equation}
for every closed interval $[a, b] \subseteq \overline{\lL}$. Moreover,
$\overline{\lL} \ominus [a, b]$ is graded of the same rank as
$\overline{\lL}$.
\end{lemma}

\noindent
\emph{Proof.} We claim that if $u < v$ holds in $\lL$ and at most one
of $u, v$ belongs to $[a, b]$, then there exists a saturated chain in
$\lL$ with minimum element $u$ and maximum element $v$ which does not
contain elements of $[a, b]$ other than $u, v$. The claim implies that
every chain in $\overline{\lL}$ which does not contain two elements of
$[a, b]$ is also a chain in $\overline{\lL} \ominus [a, b]$ or,
equivalently, that the right-hand side of (\ref{eq:lematomic}) is
contained in the left-hand side. The reverse inclusion follows from
(\ref{eq:costinclusion}). The last statement in the lemma also follows
from the claim.

To prove the claim, we apply induction on the length of the closed
interval $[u, v]$ of $\lL$. Since the result is trivial if this length
equals one, we may assume that $v$ does not cover $u$ in $\lL$. By our
assumption on $u, v$, we have $a \not\le u$ or $v \not\le b$ in $\lL$.
Suppose that $a \not\le u$. Since $\lL$ is relatively atomic, at least
two elements of the interval $[u, v]$ cover $u$. Since the meet of these
elements is $u$, at least one of them, say $w$, satisfies $a \not\le w$.
By the induction hypothesis, there exists a maximal chain
of the interval $[w, v]$ of $\lL$ which does not contain elements of
$[a, b]$ other than $v$. Adding $u$ to this chain gives the desired
maximal chain in $[u, v]$. A similar argument works if $v \not\le b$.
\qed

\bigskip
The last two parts of the following proposition are probably known;
we provide a proof for the sake of completeness. Given a simplicial
complex $\Delta$, a geometric simplicial complex $\Gamma$ realizing
$\Delta$ and a face $F \in \Delta$, we denote by $|\ast_\Delta
(F)|$ the union of the relative interiors of the geometric simplices
of $\Gamma$ which correspond to the elements of $\ast_\Delta (F)$.
We will naturally consider $|\ast_\Delta(F)|$ to be a subspace of
$|\Delta|$.

\begin{proposition} \label{prop:gorenacyclic}
Let $\lL$ be a Gorenstein* lattice. Then:
\begin{itemize}
\itemsep=0pt
\item[{\rm (i)}]
$\Delta(\overline{\lL} \ominus [a, b])$ is acyclic for every closed
interval $[a, b] \subseteq \overline{\lL}$ of positive rank.

\item[{\rm (ii)}]
$\Delta(\overline{\lL} \sm \lL_{\ge a})$ is acyclic for every $a \in
\overline{\lL}$.

\item[{\rm (iii)}]
$\Delta(\overline{\lL} \sm \lL_{\le b})$ is acyclic for every $b \in
\overline{\lL}$.
\end{itemize}
\end{proposition}

\noindent
\emph{Proof.}
To simplify notation, we set $\Delta = \Delta(\overline{\lL})$. We
first rewrite (\ref{eq:lematomic}) as
\[ \Delta(\overline{\lL} \ominus [a, b]) \ = \ \Delta \, \sm
\bigcup_{a \le x < y \le b} \ast_{\Delta} (\{x, y\}). \]
Taking geometric realizations, we get
\begin{equation} \label{eq:openstars}
|\Delta(\overline{\lL} \ominus [a, b])| \ = \ |\Delta| \, \sm
\bigcup_{a \le x < y \le b} |\ast_{\Delta} (\{x, y\})|.
\end{equation}
Since $|\Delta|$ is a homology sphere, it follows from the Lefschetz
duality theorem for homology manifolds \cite[Theorem~70.2]{Mu84} and
the long exact homology sequence of a pair that in order to prove
(i), it suffices to show that the union in the right-hand side of
(\ref{eq:openstars}) is contractible. For this, we will use a variant
of the nerve lemma \cite[Theorem~10.7]{Bj95}. Let $E$ be the set of
two-element chains $\{x, y\}$ in the closed interval $[a, b]$ and note
that $|\ast_{\Delta} (F)|$ is a nonempty open subspace of $|\Delta|$
for every $F \in E$. Note also that, for $S \subseteq E$,
\[ \bigcap_{F \in S} \ast_{\Delta} (F) \ = \ \cases{
   \varnothing, & if \ $\cup_{F \in S} \, F$ is not a chain in
   $[a, b]$, \cr
   \ast_{\Delta} (\cup_{F \in S} \, F)              & otherwise.} \]
Moreover, $|\ast_{\Delta} (G)|$ is
a cone over any point in the relative interior of $G$, and hence
contractible, for every nonempty $G \in \Delta$. Thus, by
\cite[Theorem~10.7]{Bj95}, the union in the right-hand side of
(\ref{eq:openstars}) is homotopy equivalent to the (geometric
realization of the) simplicial complex on the ground set $E$
whose faces are the sets $S \subseteq E$ for which $\cup_{F \in
S} \, F$ is a chain in $[a, b]$. Clearly, the latter complex is a
cone with apex $\{a, b\}$ and hence contractible. Therefore, so is
the union in the right-hand side of (\ref{eq:openstars}). This
completes the proof of (i).

To verify part (ii), note that $|\Delta(\overline{\lL} \sm
\lL_{\ge a})|$ is homotopy equivalent to $|\Delta| \sm |\Delta
(\lL_{\ge a})|$ by Lemma~\ref{lem:induced}. Since $\Delta(\lL_{\ge
a})$ is a cone and hence contractible, the result follows as in
the proof of part (i). Part (iii) follows from part (ii) by
passing to the dual lattice.
\qed

\bigskip
\noindent
\emph{Proof of Theorem~\ref{thm:goren}.} Let $\lL$ be a Gorenstein*
lattice and $[a, b]$ be a nonempty closed interval in $\lL$ (intervals
without subscripts in this proof are meant to be intervals in $\lL$).
By Lemma~\ref{lem:goren}, the poset $\overline{\lL} \ominus [a, b]$ is
graded of the same rank as $\overline{\lL}$.

To prove that it is Cohen-Macaulay set $\qQ := \lL \ominus [a, b]$,
consider elements $u, v \in \qQ$ with $u < v$ and note that at most one
of $u, v$ belongs to $[a, b]$. We need to show that the open interval
$(u, v)_\qQ$ has vanishing reduced homology at ranks smaller than the
rank of $(u, v)_\qQ$. This is clear if $[a, b]$ and $[u, v]$ have void
intersection, or exactly one common element, since then $[u, v]_\qQ =
[u, v]_\lL$ is a Gorenstein* lattice. Otherwise, we have $a \vee u <
b \wedge v$, and hence $a < v$ and $u < b$, in $\lL$. Note that we
may not have $a \le u$ and $v \le b$ in $\lL$. Unraveling the relevant definitions and using the claim in the proof of 
Lemma~\ref{lem:goren}, we find that
\[ (u, v)_\qQ \ = \ \cases{ (u, v) \ominus [a \vee u, b \wedge v], &
   if \ $a \not\le u$ and $v \not\le b$, \cr
   (u, v) \sm (u, b \wedge v], & if \ $a \le u$, \cr
   (u, v) \sm [a \vee u, v), & if \ $v \le b$,} \]
where the equality in the last two cases follows from the fact that the elements of the half-open intervals of $\lL$ being removed are no longer greater than $u$ (respectively, less than $v$) in $\qQ$.

Since $[u, v]$ is a Gorenstein* lattice,
Proposition~\ref{prop:gorenacyclic} implies that $(u, v)_\qQ$ is
acyclic in each case and the proof follows.
\qed

\section{Remarks}
\label{sec:rem}

Shellability of simplicial complexes and posets is a combinatorial
notion which is stronger than Cohen--Macaulayness; see, for instance,
\cite[Section~11]{Bj95}. A classical result of Bruggesser and
Mani~\cite{BM71} implies that face lattices of convex polytopes are
shellable. We call a poset $\pP$ \emph{edgewise strongly shellable}
if $\pP \ominus \iI$ is shellable of the same rank as $\pP$ for every
closed interval $\iI \subseteq \pP$.

\begin{question} \label{que:polytopes}
Are the proper parts of face lattices of convex polytopes edgewise
strongly shellable?
\end{question}

An affirmative answer has been given by the second author for Boolean
and cubical lattices.

It would be interesting to find classes of edgewise strongly
Cohen--Macaulay posets other than that provided by
Theorem~\ref{thm:goren}. Given the rich Cohen--Macaulay connectivity
properties of geometric lattices \cite{AB14} \cite{Ba82}
\cite[Section~3]{Bj80} \cite{WW86} (see the discussion in
\cite[Section~3]{Bj14}), it seems reasonable to expect that the
following question has an affirmative answer.

\begin{question} \label{que:geom}
Is the poset $\overline{\lL}$ edgewise strongly Cohen--Macaulay for
every geometric lattice $\lL$?
\end{question}

One may even guess that every 2-Cohen--Macaulay poset which forms the
proper part of a lattice is edgewise strongly Cohen--Macaulay. This
statement is false, as the following example shows. Let $\Delta$ be
the simplicial complex on the ground set $\{a, b, c, d, e\}$ consisting
of all proper subsets of $\{a, b, c, d\}$ and $\{b, c, d, e\}$ and
let $\lL$ be the face lattice of $\Delta$, meaning the poset having a
$\hat{0}$ and $\hat{1}$ for which $\overline{\lL}$ is the set of
nonempty faces of $\Delta$, ordered by inclusion. Then $\lL$ is a
lattice and $\overline{\lL}$ is 2-Cohen--Macaulay, since
$\Delta(\overline{\lL})$ is the barycentric subdivision of the
2-Cohen--Macaulay complex $\Delta$. We set $u = \{b\}$ and $v =\{b, c,
d\}$ and leave to the reader to check that the poset $(\overline{\lL}
\ominus [u, v])_{>u}$ is disconnected of rank one. This implies that
$\overline{\lL} \ominus [u, v]$ is not Cohen--Macaulay and hence that
$\overline{\lL}$ is not edgewise strongly Cohen--Macaulay.

The lattice of the previous example is
supersolvable~\cite[Example~3.14.4]{StaEC1}. Thus, the proper parts
of supersolvable lattices with nowhere zero M\"obius function are
not always edgewise strongly Cohen--Macaulay either (these posets
were shown to be doubly Cohen--Macaulay by Welker~\cite{We95}). On
the other hand, it is an interesting problem to determine the
edgewise Cohen--Macaulay connectivity of noncrossing partition
lattices and posets of injective words; see \cite{ABW07, KK13} and
references therein.

\bigskip
We close with the following question; an affirmative answer would
extend a well known property of 2-Cohen--Macaulay posets to all
edgewise 2-Cohen--Macaulay posets.

\begin{question} \label{que:mobius}
Is it true that $\widehat{\pP}$ has a nowhere zero M\"obius function
for every edgewise doubly Cohen--Macaulay poset $\pP$?
\end{question}

\end{document}